\documentclass[12pt]{article}
\usepackage[dvips]{graphicx}
\usepackage{amssymb}
\usepackage{amsmath,amsthm}
\usepackage{natbib}

\def\supp{{\rm supp}}

\usepackage{bm}
\bmdefine{\Bzero}{0}
\bmdefine{\Bone}{1}

\def\Bt{{\bf t}}
\def\Bx{{\bf x}}
\def\By{{\bf y}}
\def\Bz{{\bf z}}
\def\Ba{{\bf a}}
\def\Be{{\bf e}}

\def\Bu{{\bf u}}
\def\matrixA{{\cal A}}

\date{March, 2006}

\newtheorem{theorem}{Theorem}[section]
\newtheorem{lemma}{Lemma}[section]
\newtheorem{corollary}{Corollary}[section]

\newtheorem{definition}{Definition}[section]
\newtheorem{proposition}{Proposition}[section]

\title{Indispensable monomials of toric ideals and\\  Markov bases}
\author{
Satoshi Aoki\\
Department of Mathematics and Computer Science\\
Kagoshima University\\
Akimichi Takemura\\
Graduate School of Information Science and Technology\\
University of Tokyo\\
and\\
Ruriko Yoshida\\
Mathematics Department,  Duke University
}
\begin{document}
\maketitle
\begin{abstract}
  Extending the notion of indispensable binomials of a toric ideal
  (\cite{takemura-aoki-2004aism}, \cite{ohsugi-hibi-indispensable}),
  we define indispensable monomials of a toric ideal and establish
  some of their properties.  They are useful for
  searching indispensable binomials of a toric ideal and for proving
  the existence or non-existence of a unique minimal system of
  binomials generators of a toric ideal.  Some examples of
  indispensable monomials from statistical models for contingency tables
  are given.
\end{abstract}

\section{Introduction}
In recent years techniques of computational commutative algebra found
applications in many fields, such as optimization \cite{sturmfels1996}, 
computational biology \cite{ascb, Pachter2004a, Hosten2004}, and statistics 
\cite{Pistone2001}. Particularly,  the algebraic view of discrete statistical 
models has been applied in many statistical problems, including
conditional inference \cite{diaconis-sturmfels}, disclosure limitation 
\cite{Sullivant2005}, the maximum likelihood estimation \cite{Hosten2004}, 
and parametric inference \cite{Pachter2004a}.
{\em Algebraic statistics} is a new field, less than a decade old, and its 
term was coined by Pistone, Riccomagno and Wynn, by the title of their book 
\cite{Pistone2001}. 
Computational algebraic statistics has been very actively developed by both 
algebraists and statisticians since the pioneering work of Diaconis and 
Sturmfels \cite{diaconis-sturmfels}.  For sampling from a finite
sample space using Markov chain Monte Carlo methods, Diaconis and
Sturmfels \cite{diaconis-sturmfels} defined the notion of Markov bases and
showed that a Markov basis is equivalent to a system
of binomial generators of a toric ideal.  

In statistical applications, the number of indeterminates is often
large, and at the same time, there exists some inherent symmetry in the
toric ideal.  For this reason, we encounter computational
difficulty in applying Gr\"obner bases to
statistical problems. 
In particular, even the reduced Gr\"obner basis of a toric ideal
may contain more than several thousands elements, but one might 
be able to describe the basis concisely using symmetry.  
For example, \cite{aoki-takemura-2003metr} shows that the unique minimal
Markov bases for $3\times 3\times K, K \geq 5$ contingency tables with fixed
two-dimensional marginals contain only $6$ orbits with respect to the group
actions of permuting levels for each axis of contingency tables, while
there are $3240$, $12085$, and $34790$ elements in reduced Gr\"obner bases
for $K = 5,6, $ and $7$, respectively.  Furthermore in this example
the reduced Gr\"obner basis contains dispensable binomials and is not minimal.

Because of the difficulty mentioned above, the first two authors of this
paper have been investigating the question of minimality of
Markov bases.  In \cite{takemura-aoki-2004aism}, we defined the notion
of indispensable moves, which belong to every Markov basis.  We 
showed that there exists a unique minimal Markov basis if and only if
the set of indispensable moves forms a Markov basis.  Shortly after,
Ohsugi and Hibi investigated indispensable binomials.  They showed
that the set of indispensable binomials coincides with the
intersection of all reduced Gr\"obner basis with respect to
lexicographic term orders in \cite{ohsugi-hibi-contingency-tables-2005}.
Thus, we are interested in enumerating indispensable binomials of a given
toric ideal.  However, in general, the enumeration itself is a 
difficult problem.

This paper proposes the notion of indispensable
monomials and investigate some of their properties.  The set of
indispensable monomials contains all terms of indispensable
binomials.  Therefore if we could enumerate indispensable monomials,
then it would be straightforward to enumerate indispensable binomials.
Here it may seem that we are replacing a hard problem by a
harder one.  Computationally this may well be the case, but
we believe that the notion of indispensable monomials may be useful for
understanding indispensable binomials and finding the
existence of the unique minimal Markov basis.

In Section \ref{sec:preliminaries}  we will set notation and summarize 
relevant results from \cite{takemura-aoki-2004aism}.  In Section
\ref{sec:definition}  we will define indispensable monomials and prove 
some basic properties of the indispensable monomials.  Further
characterizations of indispensable monomials are given in Section
\ref{sec:more-properties}  and some nontrivial examples are given in
Section \ref{sec:examples}.  We will conclude with some discussions in
Section \ref{sec:discussions}.

\section{Preliminaries}
\label{sec:preliminaries}

In this section we will set appropriate notation and then
summarize main results from \cite{takemura-aoki-2004aism}.
Because of the fundamental equivalence mentioned in \cite{diaconis-sturmfels},
we use ``system of binomial generators'' and ``Markov basis'' synonymously.
Other pairs of synonyms used in this paper are (``binomial'',``move''),
(``monomial'', ``frequency vector'') and (``indeterminate'', ``cell''), as
explained below.

\subsection{Notation}
Because this paper is based on
\cite{takemura-aoki-2004aism}, we follow its notation and
terminology in statistical context.
Also we adapt some notation in
\cite{sturmfels1996, miller-sturmfels}.
Vectors, through this paper, are column vectors and $\Bx'$ denotes the
transpose of the vector $\Bx$.

Let ${\cal I}$ be a finite set of $p = |{\cal I}|$ elements. Each
element of ${\cal I}$ is called a {\em cell}.
By ordering cells, we write ${\cal I} = \{1,\ldots,p\}$ hereafter. 
A nonnegative integer
$x_{i} \in \mathbb{N} = \{0,1,\ldots\}$ is a {\em frequency} of a cell
$i$ and $\Bx = (x_1,\ldots,x_p)' \in \mathbb{N}^p$ is a {\em frequency
vector} (nonnegative integer vector). 
We write $|\Bx| = \sum_{i = 1}^{p}x_{i}$ to denote
the {\em sample size} of $\Bx$. 
In a framework of similar tests in statistical theory (see
Chapter 4 of \cite{Lehmann-tsh-3rd}), 
we consider a  $d$-dimensional sufficient statistic
defined by
\[
\Bt = \sum_{i = 1}^{p}\Ba_{i}x_{i},
\]
where $\Ba_{i} \in \mathbb{Z}^d = \{0,\pm 1,\ldots\}^d$ 
is a $d$-dimensional fixed integral column vector for 
$i = 1,\ldots,p$. 
Let $\matrixA=(\Ba_1,\ldots,\Ba_p)=(a_{ji})$ denote a $d\times p$ integral 
matrix, where $a_{ji}$ is the $j$-th element of $\Ba_{i}$. 
Then the {\em sufficient statistic} $\Bt$ is written as $\Bt = \matrixA \Bx$. 
The set of frequency vectors for a given sufficient statistic $\Bt$ is called a
{\em $\Bt$-fiber} defined by
\[
{\cal F}_{\Bt} = \{ \Bx \in \mathbb{N}^p  \mid  \Bt = \matrixA \Bx \}.
\]
A frequency vector $\Bx$ $(\in \mathbb{N}^p)$ belongs to the fiber
${\cal F}_{\matrixA\Bx}$ by definition.  We assume that a toric
ideal is homogeneous, i.e.\ there exists ${\bf w}$ such that ${\bf w}'\Ba_i=1$,
$i=1,\dots,p$.  In this case the sample size of $\Bt$ is well defined
by $|\Bt| = |\Bx|$ where $\Bx \in {\cal F}_{\Bt}$.  If the size of
${\cal F}_{\matrixA\Bx}$ is 1, i.e.\
\[
{\cal  F}_{\matrixA\Bx} = \{ \Bx \},
\]
we call $\Bx \in \mathbb{N}^p$  a {\em 1-element fiber}.
$|{\cal  F}_{\matrixA\Bx}|$ denotes the size (the
number of the elements) of the fiber ${\cal  F}_{\matrixA\Bx}$.
The {\em support} of $\Bx$
is denoted by $\supp(\Bx)= \{ i \mid x_i > 0 \}$ and the $i$-th
coordinate vector is denoted by $\Be_i=(0,\ldots,0,1,0,\ldots,0)'$,
where $1$ is in the $i$-th position.

Now, we consider the connection between contingency tables and toric 
ideals. 
Let $k[u_1,\ldots,u_p]=k[\Bu]$ denote the polynomial ring in $p$
indeterminates $u_1, \ldots, u_p$ over the field $k$.
We identify the indeterminate $u_i \in \Bu$ with the cell $i \in {\cal I}$.  
For a $p$-dimensional column vector $\Bx \in {\mathbb N}^p$  of
non-negative integers, let $\Bu^\Bx = u_1^{x_1}\cdots u_p^{x_p} \in
k[\Bu]$ denote a monomial.
For the sufficient statistic $\Bt$, we also treat 
$\Bt=(t_1,\ldots,t_d)'$ as indeterminates.
Let $k[\Bt^{\pm 1}]=k[t_1,\ldots,t_d,
t_1^{-1},\ldots,t_d^{-1}]$ denote the Laurent polynomial ring.  
Then the system of equations $\Bt = \matrixA \Bx$ is identified as the mapping
$\hat\pi: k[\Bu] \rightarrow k[\Bt^{\pm 1}]$ defined as  $x_i \mapsto
\Bt^{\Ba_i}= t_1^{a_{1i}} \cdots t_d^{a_{di}}$.  The kernel of $\hat\pi$
is denoted by $I_\matrixA={\rm ker}(\hat\pi)$ and it is the {\em toric ideal} 
associate to the matrix $\matrixA$.

For statistical applications, it is important to construct a connected
Markov chain over the given $\Bt$-fiber.
In \cite{diaconis-sturmfels}, Diaconis and Sturmfels showed that a
generator of the toric ideal $I_\matrixA$ forms a {\em Markov basis},
i.e., it can give a 
connected chain for any $\Bt$-fiber.

A $p$-dimensional integral column vector $\Bz \in \mathbb{Z}^p$ is a
{\em move} (for $\matrixA$) if
it is in the kernel of $\matrixA$, 
\[
 \matrixA\Bz = \Bzero.
\]
Let $\Bz^+ = (z_1^+,\ldots,z_p^+)'$ and
$\Bz^- = (z_1^-,\ldots,z_p^-)'$ denote the positive and negative
part of a move $\Bz$ given by
\[
z_i^+ = \max(z_i,0),\ z_i^- = -\min(z_i,0), 
\]
respectively. Then $\Bz = \Bz^+ - \Bz^-$ and $\Bz^+$ and
$\Bz^-$ are frequency vectors in the same fiber ${\cal
F}_{\matrixA\Bz^+} (=
{\cal F}_{\matrixA\Bz^-})$. Adding a move $\Bz$ to any frequency vector $\Bx$
does not change its sufficient statistic,
\[
\matrixA(\Bx + \Bz) = \matrixA\Bx,
\]
though $\Bx + \Bz$ may not necessarily be a frequency vector.
If 
adding $\Bz$ to $\Bx$ does not produce negative elements, we see that
$\Bx \in {\cal F}_{\matrixA\Bx}$ is moved to another element $\Bx +
\Bz \in {\cal F}_{\matrixA\Bx}$ by $\Bz$. In this case, we say
that a move $\Bz$ is {\em applicable to} $\Bx$. $\Bz$ is applicable to $\Bx$
if and only if $\Bx + \Bz \in {\cal F}_{\matrixA\Bx}$, and
equivalently, $\Bx \geq \Bz^-$, i.e., $\Bx - \Bz^- \in \mathbb{N}^p$.
In particular, $\Bz$ is applicable to $\Bz^-$.  We say that a move
$\Bz$ {\em contains} a frequency vector $\Bx$ if $\Bz^+=\Bx$ or $\Bz^- =
\Bx$.  The sample size of $\Bz^+$ (or $\Bz^-$) is called a {\em degree} of
$\Bz$ and denoted by
\[
\deg(\Bz) = |\Bz^+| = |\Bz^-|. 
\]
We also write $|\Bz| = \sum_{i=1}^{p}|z_i| = 2\deg(\Bz)$.

Let ${\cal B} = \{\Bz_1,\ldots,\Bz_L\}$ be a finite set of moves. Let
$\Bx$ and $\By$ be frequency vectors in the same fiber, i.e., 
$\Bx, \By \in {\cal F}_{\matrixA\Bx} (= {\cal F}_{\matrixA\By})$. 
Following \cite{takemura-aoki-2004aism}, we say that $\By$ is {\em accessible}
from $\Bx$ by ${\cal B}$ if there exists a sequence of moves
$\Bz_{i_1},\ldots,\Bz_{i_k}$ from ${\cal B}$ and $\epsilon_j \in
\{-1,+1\},\ j = 1,\ldots,k$, satisfying $\By = \Bx +
\sum_{j=1}^k\epsilon_j\Bz_{i_j}$ and $\Bx +
\sum_{j=1}^h\epsilon_j\Bz_{i_j} \in {\cal F}_{\matrixA\Bx},\ h =
1,\ldots,k-1$. The latter relation means that the move $\Bz_{i_h}$ is
applicable to $\Bx + \sum_{j=1}^{h-1}\epsilon_j\Bz_{i_j}$ for $h =
1,\ldots,k$. We write $\Bx \sim \By\ ({\rm mod}\ {\cal B})$ if $\By$ is
accessible from $\Bx$ by ${\cal B}$. An accessibility by ${\cal
B}$ is an equivalence relation in ${\cal F}_{\Bt}$ for any $\Bt$ and
each ${\cal F}_{\Bt}$ is partitioned into disjoint equivalence classes
by ${\cal B}$ (see \cite{takemura-aoki-2004aism} for detail). 
We call these equivalence classes {\em ${\cal B}$-equivalence classes} of
${\cal F}_{\Bt}$. 
Because of symmetry, we also say that $\Bx$ and $\By$
are {\em mutually accessible} by ${\cal B}$ if $\Bx \sim \By\ ({\rm mod}\
{\cal B})$. Conversely, if $\Bx$ and $\By$ are not mutually accessible
by ${\cal B}$, i.e., $\Bx$ and $\By$ are elements from two different
${\cal B}$-equivalence classes of ${\cal F}_{\matrixA\Bx}$, we say that a move
$\Bz = \Bx - \By$ {\em connects} these two equivalence classes.

A Markov basis is defined by \cite{diaconis-sturmfels} as follows. A
set of finite moves ${\cal B} = \{\Bz_1,\ldots,\Bz_L\}$ is a {\em Markov
basis} if ${\cal F}_{\Bt}$ itself forms one
${\cal B}$-equivalence class for all $\Bt$. In other words, if ${\cal B}$ is a
Markov basis, every $\Bx, \By \in {\cal F}_{\Bt}$ are mutually
accessible by ${\cal B}$ for every $\Bt$. In statistical applications,
a Markov basis makes it possible to construct a connected Markov chain
over ${\cal F}_{\matrixA\Bx}$ for any observed frequency data
$\Bx$.

Diaconis and Sturmfels \cite{diaconis-sturmfels} showed the existence of
a finite Markov basis for any $\matrixA$
and gave an algorithm to compute one. These results were obtained by
showing the fact that ${\cal B} = \{\Bz_1,\ldots,\Bz_L\}$ is a Markov
basis if and only if the set of binomials $\{\Bu^{\Bz_k^+} -
\Bu^{\Bz_k^-}, k = 1,\ldots,L\}$ is a generator of the toric ideal
$I_{\matrixA}$ associate to $\matrixA$. 
The algorithm in \cite{diaconis-sturmfels} is based on
the elimination theory of polynomial ideals and computation of a
Gr\"obner basis.

\subsection{Summary of relevant facts on indispensable moves and
  minimal Markov bases.}
In \cite{aoki-takemura-2003anz, takemura-aoki-2004aism, takemura-aoki-2005bernoulli}, we have investigated the
question on the minimality and unique minimality of Markov bases without
computing a Gr\"obner basis of $I_{\matrixA}$.
A Markov basis ${\cal B}$ is {\em minimal}
if no proper subset of ${\cal B}$ is a Markov basis. A minimal Markov
basis always exists, because from any Markov basis, we can remove
redundant elements one by one until none of the remaining elements
can be removed any further. In defining minimality 
of Markov basis, we have to be careful on signs of moves,
because minimal $\cal B$ can contain only one of $\Bz$ or $-\Bz$.
Also a minimal Markov
basis is unique if all minimal Markov bases coincide except for 
signs of their elements (\cite{takemura-aoki-2004aism}).

The structure of the unique minimal Markov
basis is given in \cite{takemura-aoki-2004aism}. Here we will summarize the 
main results of the paper without proofs.
Two particular sets of moves are important.
One is the set of moves $\Bz$ with the
same value of the sufficient statistic $\Bt = \matrixA\Bz^+ =
\matrixA\Bz^-$, namely
\[
{\cal B}_{\Bt} = \{\Bz \ |\ \matrixA\Bz^+ = \matrixA\Bz^- = \Bt\},
\]
and the other is the set of moves with degree less than or equal to $n$,
namely
\[
{\cal B}_n = \{\Bz\ |\ \deg(\Bz) \leq n\}. 
\]
Consider the ${\cal B}_{|\Bt|-1}$-equivalence classes of ${\cal
F}_{\Bt}$ for each $\Bt$. We write this
equivalence classes of ${\cal F}_{\Bt}$ as
$
{\cal F}_{\Bt} = {\cal F}_{\Bt,1} \cup \cdots \cup {\cal F}_{\Bt,K_{\Bt}}.
$

\begin{theorem}[Theorem 2.1 in \cite{takemura-aoki-2004aism}]
\label{thm:th2.1-of-takemura-aoki-2004aism}
Let ${\cal B}$ be a minimal Markov basis. 
Then for each $\Bt$, ${\cal B}\cap{\cal B}_{\Bt}$ consists of
 $K_{\Bt}-1$ moves
 connecting different ${\cal B}_{|\Bt|-1}$-equivalence classes of ${\cal
 F}_{\Bt}$, such that the equivalence classes are connected
 into a tree by these moves.

Conversely, choose any $K_{\Bt}-1$ moves
 $\Bz_{\Bt,1},\ldots,\Bz_{\Bt,K_{\Bt}-1}$ connecting different ${\cal
 B}_{|\Bt|-1}$-equivalence classes of ${\cal F}_{\Bt}$ such 
 that the equivalence classes are connected into a tree by these
 moves. Then
\[
\displaystyle\bigcup_{\Bt : K_{\Bt} \geq
 2}\{\Bz_{\Bt,1},\ldots,\Bz_{\Bt,K_{\Bt}-1}\}
\]
is a minimal Markov basis.
\end{theorem}

From  Theorem 
\ref{thm:th2.1-of-takemura-aoki-2004aism}, we immediately have
a necessarily and sufficient condition for the existence of a unique minimal
Markov basis. 

\begin{corollary}[Corollary 2.1 in \cite{takemura-aoki-2004aism}]
A minimal Markov basis is unique if and only if for each $\Bt$, ${\cal
 F}_{\Bt}$ itself forms one ${\cal B}_{|\Bt|-1}$-equivalence class or 
${\cal  F}_{\Bt}$ is a two-element fiber.
\end{corollary}

This condition is explicitly expressed by {\em indispensable moves}.

\begin{definition}
A move $\Bz = \Bz^+ - \Bz^-$ is called indispensable if 
$\Bz^+$ and $\Bz^-$ form a two-element fiber, i.e., 
the fiber ${\cal F}_{\matrixA\Bz^+} (= {\cal F}_{\matrixA\Bz^-})$ is
 written as ${\cal F}_{\matrixA\Bz^+} = \{\Bz^+,\Bz^-\}$.
\end{definition}

From the above definition and the structure of a minimal Markov basis,
one can show that every indispensable move belongs to each Markov basis
(Lemma 2.3 in \cite{takemura-aoki-2004aism}). Furthermore, by the
correspondence between moves and binomials, we define an 
{\em indispensable binomial}.

\begin{definition}
A binomial $\Bu^{\Bz} = \Bu^{\Bz^+} - \Bu^{\Bz^-}$ is indispensable if
 every system
 of binomial generators of $I_{\matrixA}$ contains $\Bu^{\Bz}$ or $-\Bu^{\Bz}$. 
\end{definition}

Clearly, a binomial $\Bu^{\Bz}$ is indispensable if and only if a move
$\Bz$ is indispensable.
By definition, a set of indispensable moves plays an important role to 
determine the uniqueness of a minimal Markov basis:

\begin{lemma}[Corollary 2.2 in \cite{takemura-aoki-2004aism}]
\label{lemmaMB}
The unique minimal Markov basis exists if and only if the set of
 indispensable moves forms a Markov basis. In this case, the set of
 indispensable moves is the unique minimal Markov basis.
\end{lemma}

Ohsugi and Hibi further investigated indispensable moves 
\cite{ohsugi-hibi-contingency-tables-2005, ohsugi-hibi-indispensable}.

\begin{theorem}[Theorem 2.4 in \cite{ohsugi-hibi-contingency-tables-2005}]
A binomial $\Bu^{\Bz}$ is indispensable if and only if either 
$\Bu^{\Bz}$ or $-\Bu^{\Bz}$ belongs to the reduced Gr\"obner basis of
 $I_{\matrixA}$ for any lexicographic term order on $k[\Bu]$.
\end{theorem}

One can find more details in \cite{ohsugi-hibi-indispensable}.

\section{Definition and some basic properties of indispensable monomials}
\label{sec:definition}

In this section we will define indispensable monomials.  
Then we will show two other equivalent conditions for a monomial to be
indispensable. 
We will also prove analogous to
Theorem 2.4 in \cite{ohsugi-hibi-contingency-tables-2005}, that the set of
indispensable monomials is characterized as the intersection of monomials
in reduced Gr\"obner bases with respect to all lexicographic
term orders. 
Hereafter, we say that a  Markov basis $\cal B$ contains $\Bx$ if
it contains a move $\Bz$ containing $\Bx$ by  abusing the
terminology. 

Firstly we will define an {\em indispensable monomial}.

\begin{definition}
\label{def:1}
A  monomial $\Bu^\Bx$ is {\it indispensable} if
every system of binomial generators of $I_{\matrixA}$ contains a binomial
$f$  such that $\Bu^\Bx$ is a term of $f$.
\end{definition}

From this definition, any Markov basis contains all indispensable
monomials. Therefore the set of indispensable monomials is finite.
Note that both terms of an indispensable binomial $\Bu^{\Bz^+} -
\Bu^{\Bz^-}$ are indispensable monomials, but the converse does not
hold in general.

Now we will present an alternative definition.

\begin{definition}
$\Bx$ is a minimal multi-element if $|{\cal  F}_{\matrixA\Bx}|\ge 2$ and 
$|{\cal  F}_{\matrixA(\Bx-\Be_i)}|=1$ for every $i\in \supp(\Bx)$.
\end{definition}

\begin{theorem}
\label{thm:1}
$\Bx$ is an indispensable monomial if and only if $\Bx$ is a minimal
multi-element.
\end{theorem}

\begin{proof}
First, we suppose that $\Bx$ is a minimal multi-element and want to show that
it is an indispensable monomial. Let $n=|\Bx|$  and 
consider the fiber ${\cal F}_{\matrixA\Bx}$.  
We claim that $\{\Bx \}$ forms a single ${\cal B}_{n-1}$-equivalence
class.  In order to show this,  we argue by contradiction.  If 
$\{\Bx \}$ does not form a single ${\cal B}_{n-1}$-equivalence class, then
there exists a  move $\Bz$ with degree less than or equal to $n-1$, such that
\[
\Bx + \Bz = (\Bx - \Bz^-) + \Bz^+ \in {\cal F}_{\matrixA\Bx}.
\]
Since $|\Bx|=n$,  $|\Bz^-|\le n-1$, we have  
$\Bzero \neq \Bx - \Bz^-$ and 
\[
\supp (\Bx) \cap \supp(\Bx + \Bz) \neq \emptyset.
\]
Therefore we can choose  $i\in \supp (\Bx) \cap \supp(\Bx + \Bz)$ such that
\[
\matrixA (\Bx - \Be_i) = \matrixA (\Bx + \Bz - \Be_i), \quad \Bx -
 \Be_i \neq \Bx + \Bz - \Be_i.
\]
This shows that $|{\cal  F}_{\matrixA(\Bx-\Be_i)}|\ge 2$, which contradicts the
assumption that $\Bx$ is a minimal multi-element.
Therefore we have shown that $\{\Bx \}$ forms a single ${\cal
  B}_{n-1}$-equivalence class.
Since we are assuming that $|{\cal F}_{\matrixA\Bx}|\ge 2$, there exists
some other ${\cal B}_{n-1}$-class in ${\cal F}_{\matrixA\Bx}$.  By
Theorem \ref{thm:th2.1-of-takemura-aoki-2004aism}, 
each Markov basis has to contain a move connecting a one element
 equivalence class $\{\Bx\}$ to other equivalence classes of ${\cal
 F}_{\matrixA\Bx}$, which implies that each Markov basis has to 
contain a move $\Bz$ containing $\Bx$.
We now have shown that each minimal multi-element has to be contained in each
Markov basis, i.e., a minimal multi-element is an indispensable monomial.

Now we will show the converse. It suffices to show that if $\Bx$ is not a
minimal multi-element, then $\Bx$ is a dispensable monomial.
Suppose that $\Bx$ is not a minimal multi-element.  If $\Bx$ is a
1-element ($|{\cal F}_{\matrixA\Bx}|=1$), obviously it is 
dispensable.  Hence assume $|{\cal F}_{\matrixA\Bx}|\ge 2$.  In the case
that ${\cal F}_{\matrixA\Bx}$ is a single ${\cal B}_{n-1}$-equivalence
class, no move containing $\Bx$ is needed in a minimal Markov basis by
 Theorem \ref{thm:th2.1-of-takemura-aoki-2004aism}. 
Therefore we only need to consider the case that ${\cal F}_{\matrixA\Bx}$
contains more than one ${\cal B}_{n-1}$-equivalence classes.  Because
$\Bx$ is not a minimal multi-element, there exists some $i\in \supp(\Bx)$
such that $|{\cal F}_{\matrixA(\Bx-\Be_i)}|\ge 2$. 
Then there exists $\By \neq \Bx-\Be_i$, such that $\matrixA\By
 =\matrixA(\Bx-\Be_i)$.  
Noting that $|\By|=|\Bx-\Be_i|=n-1$, a move of the form
\[
\Bz = \By -  (\Bx-\Be_i)  = (\By + \Be_i) - \Bx
\]
satisfies $0<{\rm deg}(\Bz) \le n-1$. Then 
\[
\By + \Be_i = \Bx + \Bz
\]
and $\Bx$ and  $\By+\Be_i$ belong to the same
${\cal   B}_{n-1}$-equivalence class of ${\cal F}_{\matrixA\Bx}$. 
Since $\Bx \neq \By + \Be_i$, Theorem
 \ref{thm:th2.1-of-takemura-aoki-2004aism} states that 
we can construct a minimal Markov
basis containing $\By+\Be_i$, but not containing $\Bx$.
Therefore  $\Bx$ is a dispensable monomial.
\end{proof}

We will give yet another definition.

\begin{definition}
\label{def:3}
  $\Bx$ is a minimal $i$-lacking 1-element if
$|{\cal F}_{\matrixA\Bx}|=1$, 
 $|{\cal F}_{\matrixA(\Bx+\Be_i)}|\ge 2$ and 
$|{\cal F}_{\matrixA(\Bx+\Be_i- \Be_j)}|=1$ for each $j\in \supp(\Bx)$.
\end{definition}

We then have the following result.

\begin{theorem}
\label{thm:2}  
The following three conditions are equivalent
1) $\Bx$ is an indispensable monomial, 2) for each
$i\in\supp(\Bx)$, $\Bx-\Be_i$ is a minimal $i$-lacking 1-element, 3) for
some $i\in\supp(\Bx)$, $\Bx-\Be_i$ is a minimal $i$-lacking 1-element.
\end{theorem}

By the previous theorem we can replace the condition 1) by the
condition that $\Bx$ is a minimal multi-element.

\begin{proof}
$1) \Rightarrow 2)$.
Suppose that $\Bx$ is a minimal multi-element.  Then for any $i\in
\supp(\Bx)$, $\Bx-\Be_i$ is a 1-element and 
$|{\cal F}_{\matrixA((\Bx-\Be_i)+\Be_i)}|=
|{\cal F}_{\matrixA\Bx}|\ge 2$.  If $\Bx-\Be_i$ is not a minimal
$i$-lacking 1-element, then for some $j\in \supp(\Bx-\Be_i)$, 
$|{\cal F}_{\matrixA(\Bx- \Be_j)}| \ge 2$.  However 
$j\in \supp(\Bx-\Be_i) \subset \supp(\Bx)$ and 
$|{\cal F}_{\matrixA(\Bx- \Be_j)}| \ge 2$ contradicts the assumption that 
$\Bx$ is a minimal multi-element.
It is obvious that $2) \Rightarrow 3)$.  

Finally we will prove $3) \Rightarrow 1)$.
Suppose that for some $i\in \supp(\Bx)$, $\Bx- \Be_i$ is a
minimal $i$-lacking 1-element.  Note
that $|{\cal F}_{\matrixA \Bx}|=|{\cal
  F}_{\matrixA((\Bx-\Be_i)+\Be_i)}|\ge 2$.  Now consider $j\in
\supp(\Bx)$. If $j\in \supp(\Bx-\Be_i)$ then $|{\cal
  F}_{\matrixA(\Bx-\Be_j)}|=|{\cal F}_{\matrixA((\Bx-\Be_i) + \Be_i -
  \Be_j)}|=1$.  On the other hand if $j\not\in \supp(\Bx-\Be_i)$, then
$j=i$ because $j\in \supp(\Bx)$.  In this
case $|{\cal F}_{\matrixA(\Bx-\Be_i)}|=1$.  This shows that $\Bx$ is a
minimal multi-element.
\end{proof}

Theorem \ref{thm:2} suggests the following:  Find any 1-element $\Bx$.  
It is often the case
that each $\Be_i$, $i=1,\ldots,p$, is a 1-element. Randomly choose $1\le
i\le p$ and check whether $\Bx+\Be_i$ remains to be a 1-element.  Once
$|{\cal F}_{\Bx+\Be_i}| \ge 2$, then subtract $\Be_j$'s, $j\neq i$, one by
one from $\Bx$ such that it becomes a minimal $i$-lacking 1-element. 
We can apply this procedure to  finding indispensable
monomials of some actual statistical problem.

For the rest of this section we will illustrate this procedure with an
example of a $2\times 2\times 2$ contingency table.
Consider the following problem where $p = 8, d = 4$ and $\matrixA$ is given
 as
\[
\matrixA = \left(
\begin{array}{cccccccc}
1 & 1 & 1 & 1 & 1 & 1 & 1 & 1\\
1 & 1 & 1 & 1 & 0 & 0 & 0 & 0\\
1 & 1 & 0 & 0 & 1 & 1 & 0 & 0\\
1 & 0 & 1 & 0 & 1 & 0 & 1 & 0
\end{array} 
\right)\ .
\]
In statistics this is known as the complete independence model of $2\times
 2\times 2$ contingency tables. To see the direct product structure of
 ${\cal I}$ explicitly, we write indeterminates as
\[
 \Bu = (u_{111},u_{112},u_{121},u_{122},u_{211},u_{212},u_{221},u_{222}).
\]
To find indispensable monomials for this problem, we start with the
 monomial $\Bu^{\Bx} = u_{111}$ and consider $\Bx + \Be_i, i \in {\cal
 I}$. Then we see that
\begin{itemize}
\item $u_{111}^2, u_{111}u_{112}, u_{111}u_{121}, u_{111}u_{211}$ are
      $1$-elements,
\item $u_{111}u_{122}, u_{111}u_{212}, u_{111}u_{221}$ are
      $2$-elements and
\item $u_{111}u_{222}$ is a $4$-element.
\end{itemize}
From this, we found four indispensable monomials,
$u_{111}u_{122}, u_{111}u_{212}, u_{111}u_{221}$ and $u_{111}u_{222}$,
since each of $u_{122}, u_{212}, u_{221}, u_{222}$ is a $1$-element. 

Starting from the other monomials, similarly, we can find the following list
 of indispensable monomials,
\begin{itemize}
\item $u_{111}u_{122}, u_{111}u_{212}, u_{111}u_{221}, u_{112}u_{121},
 u_{112}u_{211}, u_{112}u_{222},
 u_{121}u_{222}, u_{121}u_{211},\\
 u_{122}u_{221}, u_{122}u_{212}, u_{211}u_{222}, u_{212}u_{221}$, each
      of which is a $2$-element monomial, and
\item $u_{111}u_{222}, u_{112}u_{221},u_{121}u_{212},u_{122}u_{211}$,
      each of which is a $4$-element monomial.
\end{itemize}

The next step is to consider the newly produced $1$-element monomials, \\
$u_{111}^2, u_{111}u_{112}, u_{111}u_{121}, u_{111}u_{211}$ and so on. 
For each of these monomials, consider adding $\Be_i, i \in {\cal I}$
 one by one, checking whether they are multi-element or not. For example, we
 see that the monomials such as
\[
 u_{111}^3, u_{111}^2u_{112}, u_{111}u_{112}^2,\ldots
\]
are again $1$-element monomials (and we have to consider these
 $1$-element monomials in the
 next step). On the other hand, monomials such as
\[
 u_{111}^2u_{122},  u_{111}u_{112}u_{122}, u_{111}^2u_{222},
 u_{111}u_{112}u_{221},\ldots
\]
are multi-element monomials. However, it is seen that they are
 not minimal multi-element, since
\[
 u_{111}u_{122},  u_{112}u_{122}, u_{111}u_{222}, u_{112}u_{221},\ldots
\]
are not $1$-element monomials.

To find all indispensable monomials for this problem, we have to 
 repeat the above procedure for monomials of degree $4,5,\ldots$.  Note 
 that this procedure never stops since there are
 infinite $1$-element monomials, such as
\[
 u_{111}^n,\ u_{111}^nu_{112}^m,\ldots
\]
for arbitrary $n,m$. 
This is the same difficulty mentioned in Section 2.2 in
\cite{takemura-aoki-2004aism}. 
Since indispensable monomials belong to any Markov basis, in
particular to the Graver basis, Theorem 4.7 in \cite{sturmfels1996}
gives an upper bound for the degree of indispensable monomials
and we can stop at this bound.

Finally we will prove the following theorem, which is analogous to Theorem
2.4 in \cite{ohsugi-hibi-contingency-tables-2005} but much easier to prove,
since it focuses on a single monomial (rather than a binomial).
We need to reproduce only a part of the proof for Theorem
2.4 in \cite{ohsugi-hibi-contingency-tables-2005}.

\begin{theorem}
A monomial $\Bx$ is indispensable if 
for every lexicographic 
order $<_{\rm lex}$ the reduced Gr\"obner basis with respect to $<_{\rm lex}$
contains $\Bx$.
\end{theorem}

\begin{proof} It suffices to show that if a monomial $\Bx$ is 
  dispensable, then there exists a lexicographic term order $<_{\rm
    lex}$ such that the reduced Gr\"obner basis ${\cal B}_{<_{\rm
      lex}}$ does not contain $\Bx$.  Note that the positive part and
  negative part of a move belong to the same fiber.  Therefore if
  $\Bx$ is a 1-element, then no Markov basis contains $\Bx$. In
  particular no Gr\"obner basis contains $\Bx$.  Therefore we only
  need to consider $\Bx$ such that $|{\cal F}_{\matrixA\Bx}|\ge 2$.

  Since we are assuming that $\Bx$ is dispensable, there exists a
  Markov basis ${\cal B}$, which does not contain $\Bx$.  Then there
  exists a move $\Bz=\Bz^+ - \Bz^-\in {\cal B}$, with $\Bz$'s sign
  changed if necessary, such that $\Bz$ is applicable to $\Bx$, i.e., 
  $\Bx \ge \Bz^-$.  Since $\cal B$ does not
  contain $\Bx$, $\Bz^- \neq \Bx$ and hence $\Bz^-$ is strictly
  smaller than $\Bx$.  Now choose $<_{\rm lex}$ such that the initial
  term of $\Bz$ is $\Bz^-$.  Then the reduced Gr\"obner basis with
  respect to $<_{\rm lex}$ does not contain $\Bx$.
\end{proof}

\section{Further properties of indispensable monomials}
\label{sec:more-properties}

In the previous section we gave some basic characterizations of
indispensable monomials.  In this section we will show further properties of
indispensable moves in terms of minimal Markov
bases in \cite{takemura-aoki-2004aism} and a {\em norm-reducing
Markov basis} in \cite{takemura-aoki-2005bernoulli}.

Firstly, we will state the following lemma, which is already
implicitly used in the proof of Theorem \ref{thm:1}.

\begin{lemma}
\label{lem:basic}
A monomial $\Bu^\Bx$ is indispensable if and only if ${\cal
  F}_{\matrixA \Bx}$ contains more than one ${\cal
  B}_{|\Bx|-1}$-equivalence class and the one-element set $\{\Bx\}$
forms a ${\cal B}_{|\Bx|-1}$-equivalence class.
\end{lemma}

\begin{proof}
Suppose that 
${\cal  F}_{\matrixA \Bx}$ contains more than one ${\cal
  B}_{|\Bx|-1}$-equivalence class and the one-element
set $\{\Bx\}$ forms a ${\cal B}_{|\Bx|-1}$-equivalence class.
Then by Theorem 2.1 in \cite{takemura-aoki-2004aism},  
every Markov basis has to connect $\Bx$ with some other 
${\cal B}_{|\Bx|-1}$-equivalence class of ${\cal F}_{\matrixA \Bx}$.
Therefore $\Bx$ has to appear as a positive part or a negative part of
some move $\Bz$ of the Markov basis.

Conversely, we show that if ${\cal
  F}_{\matrixA \Bx}$ contains just one  ${\cal
  B}_{|\Bx|-1}$-equivalence class or the equivalence class containing
$\Bx$ contains some other vector $\By$, then $\Bx$ is 
dispensable. In the former case, ${\cal   F}_{\matrixA \Bx}$ is
already connected by moves of degree less than or equal to $|\Bx|-1$
and no minimal Markov basis contains a move having $\Bx$ as the
positive or the negative part.  On the other hand if $\By \neq \Bx$
belongs to the same ${\cal  B}_{|\Bx|-1}$-equivalence class, then 
by Theorem 2.1 in \cite{takemura-aoki-2004aism},  
there exists a minimal Markov basis
involving $\By$ and not $\Bx$.  Therefore $\Bx$ is dispensable.
\end{proof}

From \cite{takemura-aoki-2004aism} it follows that the moves of all
minimal Markov bases belong to a common set of fibers. Also,
we defined the minimum fiber Markov bases ${\cal B}_{\rm MF}$ 
in \cite{takemura-aoki-2005bernoulli} as
\[
 {\cal B}_{\rm MF} = \{\Bz = \Bz^+ - \Bz^-\ |\ \Bz^+ \not\sim \Bz^-\
 ({\rm mod}\ {\cal B}_{|\Bz|-1})\}.
\]
Based on Lemma \ref{lem:basic}
now we will prove four propositions concerning the fibers in ${\cal B}_{\rm
  MF}$. In the following four propositions, an equivalence class of a
fiber ${\cal F}_{\Bt}$ means a ${\cal B}_{|\Bt|-1}$-equivalence
class of ${\cal F}_{\Bt}$.

\begin{proposition}
\label{prop:1}
The following three conditions are equivalent: 1) all equivalence classes
of all fibers of ${\cal   B}_{\rm MF}$  are singletons, 2)
there exists a minimal Markov basis, such that all  monomials contained
in the basis are indispensable.
3) for all  minimal Markov bases,  all monomials contained in the
basis are indispensable. 
\end{proposition}

\begin{proof} Obviously $3) \Rightarrow 2)$.  $2) \Rightarrow 1)$
  follows from Lemma \ref{lem:basic} because a minimal basis has to
  connect all equivalence classes of each fiber of ${\cal B}_{\rm MF}$
  into a tree.  To show that $1) \Rightarrow 3)$, we again use the
  fact that a minimal basis has to connect all equivalence classes of
  each fiber of ${\cal B}_{\rm MF}$ into a tree. If all 
  equivalence classes of a fiber are singletons, then both terms of a
  move connecting two equivalence classes are indispensable.  This
  completes the proof.
\end{proof}

\begin{proposition}
\label{prop:2}
There exists a minimal Markov basis such that each move of the basis
contains an indispensable monomial if and only if
each fiber of ${\cal B}_{\rm MF}$ contains a singleton equivalence class.
\end{proposition}

\begin{proof} Let $\cal B$
  be a minimal Markov basis such that each move of $\cal B$ contains
  an indispensable monomial.  This monomial forms a singleton
  equivalence class.  Therefore each fiber of ${\cal B}_{\rm MF}$
  contains a singleton equivalence class. Conversely if each fiber of
  ${\cal B}_{\rm MF}$ contains a singleton equivalence class, we can
  construct a tree which connects each equivalence class of the fiber to the
  singleton equivalence class.  Then the resulting minimal Markov
  basis has the property that each move of the basis contains an
  indispensable monomial.
\end{proof}

\begin{proposition}
\label{prop:3}
Every move of any minimal Markov basis contains an indispensable
monomial if and only if all but one equivalence classes of each fiber
of ${\cal   B}_{\rm MF}$ are singletons.
\end{proposition}

\begin{proof}
If all but one equivalence classes of each fiber
of ${\cal   B}_{\rm MF}$ are singletons, then in connecting these
equivalence classes into a tree, each move has to contain an
indispensable monomial.
On the other hand if there exist two non-singleton equivalence classes
in a fiber, then we can construct a minimal Markov basis containing a
move connecting these two equivalence classes. This move does not
contain an indispensable monomial.
\end{proof}

Next, we consider indispensable monomials in terms of
norm-reduction introduced in \cite{takemura-aoki-2005bernoulli}.
We will give a definition of a norm-reducing Markov basis here (see
\cite{takemura-aoki-2005bernoulli} for detail).
\begin{definition}
A set of moves ${\cal B}$ is 1-norm reducing if for all
 $\Bt$ and for all
 $\Bx, \By \in {\cal F}_{\Bt}$ with $\Bx \neq \By$, there exist some $\Bz \in
 {\cal B}$ and $\epsilon \in \{-1,+1\}$  satisfying either 
\[
|(\Bx + \epsilon\Bz) - \By| < |\Bx - \By|
\]
or
\[
|\Bx - (\By + \epsilon\Bz)| < |\Bx - \By|.
\]
\end{definition}
It is easy to show that, if ${\cal B}$ is 1-norm reducing, then it is a
 Markov basis (see Proposition 1 in \cite{takemura-aoki-2005bernoulli}).
Therefore we call ${\cal B}$ a {\em 1-norm reducing Markov basis} if it is
 1-norm reducing. An example of 1-norm reducing Markov basis is the
 Graver basis (see Proposition 2 in \cite{takemura-aoki-2005bernoulli}).
Now we will give a characterization of indispensable monomials in terms of 
the norm reduction. 

\begin{proposition}
\label{prop:4}
A move, whose both monomials are indispensable,
belongs to each 1-norm reducing Markov basis.
\end{proposition}

\begin{proof} Let $\Bz = \Bz^+ - \Bz^-$ be a move such that both
  $\Bz^+$ and $\Bz^-$ are indispensable, i.e., $\{\Bz^+\}$ and
  $\{\Bz^-\}$ are singleton equivalence classes of a fiber.  If a
  Markov basis ${\cal B}$ does not contain $\Bz$, then it can not
  decrease the distance between $\Bz^+$ and $\Bz^-$.  Therefore ${\cal
    B}$ is not 1-norm-reducing.
\end{proof}

Finally, we give a further definition, which is similar to a minimal
multi-element.
\begin{definition}
${\cal F}_{\Bt}$ is a minimal multi-element fiber if 
\{
$|{\cal F}_{\Bt}| \geq 2$ and $|\Bt| = 1$
\}
or 
\{
$|{\cal F}_{\Bt}| \geq 2, |{\cal F}_{\Bt_1}| = |{\cal F}_{\Bt_2}| = 1$
for any $\Bt = \Bt_1 + \Bt_2$ satisfying $|{\cal F}_{\Bt_1}|, |{\cal
 F}_{\Bt_2}| \geq 1$
\}.
\end{definition}
The meaning of this definition is as follows. Suppose $|{\cal
F}_{\Bt_1}| \geq 2$ and $\Bx_1, \Bx_2 \in {\cal F}_{\Bt_1}$. Then for
any $\Bt_2 = \matrixA\Bx_3$, $|{\cal F}_{\Bt_1 + \Bt_2}| \geq 2$ follows since $\Bx_1 +
\Bx_3, \Bx_2 + \Bx_3 \in {\cal F}_{\Bt_1 + \Bt_2}$. 
Note that the former case, $|{\cal F}| \geq 2, |\Bt| = 1$, corresponds
to the case ${\cal F}_{\Bt} \ni \Be_i, \Be_j,\ldots$, for some
$i,j,\ldots$. One of the situations that this special case appears is
$\matrixA = [1,\ldots,1]$. Hereafter, we only consider the problem that 
$|\Bt| \geq 2$ holds for every minimal multi-element fiber ${\cal
F}_{\matrixA\Bt}$. In other words, we assume that every $\Be_i$ is a
one-element. In this case, minimal multi-element fiber is characterized
as follows.
\begin{proposition}
${\cal F}_{\Bt}$ is a minimal multi-element fiber if and only if all the
 elements in ${\cal F}_{\Bt}$ are indispensable monomials.
\end{proposition}

\begin{proof} 
Suppose all the elements in ${\cal F}_{\Bt}$ are indispensable monomials
 and $\Bt = \Bt_1 + \Bt_2$ where $|{\cal F}_{\Bt_1}| \geq 2, |{\cal F}_{\Bt_2}|
 \geq 1$. Write $\Bx_1, \Bx_2 \in {\cal F}_{\Bt_1}$ and $\Bx_3 \in {\cal
 F}_{\Bt_2}$. In this case, both $\Bx_1 + \Bx_3$ and $\Bx_2 + \Bx_3$
 are in ${\cal F}_{\Bt}$ and therefore indispensable monomials by the
 assumption. However, for any $i \in {\rm supp}(\Bx_3)$, 
$\Bx_1 + \Bx_3 - \Be_i$ and $\Bx_2 + \Bx_3 - \Be_i$ are again in the
 same fiber, which contradicts the assumption that $\Bx_1 + \Bx_3$ and
 $\Bx_2 + \Bx_3$ are minimal multi-elements. 

Conversely, suppose ${\cal F}_{\Bt}$ is a minimal multi-element fiber
 and $\Bx \in {\cal F}_{\Bt}$ is dispensable. In this case, since $\Bx$
 is not a minimal multi-element, there exists some $\Be_i$ satisfying
 $|{\cal F}_{\matrixA(\Bx - \Be_i)}| \geq 2$. Therefore we have
 $\matrixA\Bx = \matrixA(\Bx - \Be_i) + \matrixA\Be_i$, i.e., $\Bt =
 \Bt_1 + \Bt_2$ where $\Bt_1 = \matrixA(\Bx - \Be_i)$ and $\Bt_2 =
 \matrixA\Be_i$, which contradicts the assumption that ${\cal F}_{\Bt}$
 is a minimal multi-element fiber.

\end{proof}

\section{Examples}
\label{sec:examples}
In this section, we will give some indispensable monomials and
dispensable monomials in minimal Markov bases for some 
statistical models. 
As is stated in \cite{takemura-aoki-2004aism}, there are some models
where a minimal Markov basis is uniquely determined, and some models
where it is not uniquely determined. Furthermore, by considering the
indispensability of monomials contained in minimal Markov bases, we
can classify Markov bases by the indispensability of monomials as follows.
\begin{itemize}
\item Case 1.\ A minimal Markov basis is uniquely determined, i.e., the
      set of indispensable moves forms a Markov basis.
\item Case 2.\ A minimal Markov basis is not uniquely determined, but all
      the monomials in minimal Markov bases are the same and
      indispensable. In this case, all  equivalence classes of each
      fiber of ${\cal B}_{\rm MF}$ are singletons.
\item Case 3.\ A minimal Markov basis is not uniquely determined, and
      they contain some moves where their positive or negative parts are
      dispensable monomials. In this case, some 
      equivalence classes of some fiber of ${\cal B}_{\rm MF}$ are not
      singletons.
\end{itemize}
We will show examples for Case 2 and Case 3 in this section.
As for Case 1, the set of the positive
and negative parts of indispensable binomials is the set of
indispensable monomials. 
One of the most simple examples for Case 1 is an 
independence model of two-way contingency tables. A quite difficult
example 
is a no three-factor interaction model of three-way contingency tables,
i.e., the case that $\matrixA\Bx$ is the two-dimensional marginal totals
of three-way contingency tables $\Bx$. For this example, minimal Markov
bases for some small sizes of $\Bx$ is shown to be unique (see
\cite{aoki-takemura-2003anz} for example of $3\times 3\times K$ case). 
Indispensable monomials for Case 1 clearly coincide the 
positive and  negative parts of indispensable binomials.

\subsection{Examples of Case 2}
\paragraph*{One-way contingency tables with fixed totals.}
\ First we consider the simplest example given by $\matrixA = 1_{p}', p > 2$, 
where
$1_{p} = (1,\ldots,1)'$ is the $p$ dimensional vector consisting $1$'s. 
As is shown in \cite{takemura-aoki-2004aism}, minimal Markov bases for
this problem contain dispensable moves only, which connect $p$ elements,
\[
 \{u_{1},u_{2},\ldots,u_{p}\}
\]
into a tree. It is also obvious that these $p$ monomials are all
indispensable.

\paragraph*{Complete independence models of three-way contingency
    tables.} \ We will show a generalization of the problem considered 
at the end of Section \ref{sec:definition}.

\begin{figure}[htbp]
\begin{center}
\includegraphics[width=50mm]{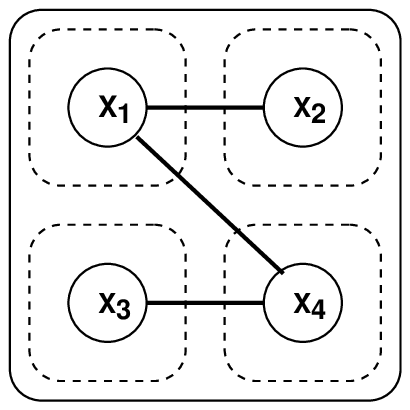}
\caption{Illustration of the 4-element fiber ${\cal
F}_{i_1i_2j_1j_2k_1k_2}$ of the complete independence model of
 three-way contingency tables. The four circles are four monomials, where
 $\Bx_1 = u_{i_1j_1k_1}u_{i_2j_2k_2}, \Bx_2 = u_{i_1j_1k_2}u_{i_2j_2k_1},
 \Bx_3 = u_{i_1j_2k_1}u_{i_2j_1k_2},$ and 
$ \Bx_4 = u_{i_1j_2k_2}u_{i_2j_1k_1}$. 
Each monomial forms 
 ${\cal B}_{1}$-equivalence class of the fiber by itself represented by the 
dotted square. We will use this convention in all forthcoming figures.  
The thick
 lines mean a choice of three dispensable moves, $\{\Bx_1-\Bx_2, \Bx_1-\Bx_4,
 \Bx_3-\Bx_4\}$, which is an example of choices for constructing a minimal
 Markov basis.}
\label{fig:4-element-monomials}
\end{center}
\end{figure}

Let $\Bx$ be a
frequency vector for $I\times J\times K$ contingency tables and let
\[
 {\cal I} = \{ ijk \ |\ 1\leq i\leq I, 1 \leq j\leq J, 1\leq k\leq K\}.
\]
$\matrixA$ is given as
\[
 \matrixA = \left[
\begin{array}{c}
1_{I}' \otimes 1_{J}' \otimes E_{K}\\
1_{I}' \otimes E_{J} \otimes 1_{K}'\\
E_{K} \otimes 1_{J}' \otimes 1_{K}'
\end{array}
\right],
\]
where 
$E_n$ is the $n\times n$ identity
matrix. The minimum fiber Markov basis for this problem is given in 
\cite{takemura-aoki-2005bernoulli} as

\[\begin{array}{l}
{\cal B}_{\rm{MF}} = {\cal B}_{\rm{IDP}} \cup {\cal B}^{*},\\
{\cal B}_{\rm{IDP}} = \{u_{ij_1k_1}u_{ij_2k_2} - u_{ij_1k_2}u_{ij_2k_1},\
j_1\neq j_2,\ k_1\neq k_2\}\\
\hspace*{15mm}\cup\ \{u_{i_1jk_1}u_{i_2jk_2} - u_{i_1jk_2}u_{i_2jk_1},\
i_1\neq i_2,\ k_1\neq k_2\}\\
\hspace*{15mm}\cup\ \{u_{i_1j_1k}u_{i_2j_2k} - u_{i_1j_2k}u_{i_2j_1k},\
i_1\neq i_2,\ j_1\neq j_2\},\\
{\cal B}^{*} = \{u_{i_1j_1k_1}u_{i_2j_2k_2} -
u_{i_1j_1k_2}u_{i_2j_2k_1},\
u_{i_1j_1k_1}u_{i_2j_2k_2} - u_{i_1j_2k_1}u_{i_2j_1k_2},\\
\hspace*{15mm} u_{i_1j_1k_1}u_{i_2j_2k_2} -
u_{i_1j_2k_2}u_{i_2j_1k_1},\
u_{i_1j_1k_2}u_{i_2j_2k_1} - u_{i_1j_2k_1}u_{i_2j_1k_2},\\
\hspace*{15mm} u_{i_1j_1k_2}u_{i_2j_2k_1} -
u_{i_1j_2k_2}u_{i_2j_1k_1},\
u_{i_1j_2k_1}u_{i_2j_1k_2} - u_{i_1j_2k_2}u_{i_2j_1k_1},\\
\hspace*{75mm} i_1\neq i_2, j_1\neq j_2, k_1\neq k_2\}.
\end{array}
\]
Here, ${\cal B}_{\rm{IDP}}$ is the set of indispensable moves. ${\cal
  B}^{*}$ is the set of all degree 2 moves which connect all 
  elements of the four-elements fiber
\[\begin{array}{ccl}
{\cal F}_{i_1i_2j_1j_2k_1k_2}
 & = & \{\Bx = \{x_{ijk}\} \ |\ x_{i_1\cdot\cdot} =
x_{i_2\cdot\cdot} =
x_{\cdot j_1\cdot} = x_{\cdot j_2\cdot} =
x_{\cdot\cdot k_1} = x_{\cdot\cdot k_2} = 1\}\\
& = & \{u_{i_1j_1k_1}u_{i_2j_2k_2},\ u_{i_1j_1k_2}u_{i_2j_2k_1},\
u_{i_1j_2k_1}u_{i_2j_1k_2},\ u_{i_1j_2k_2}u_{i_2j_1k_1}\}.
\end{array}
\]
The minimal Markov basis in this case consists of ${\cal B}_{\rm{IDP}}$ and
three moves for each $i_1\neq i_2,j_1\neq j_2, $ and $k_1\neq k_2$, 
which connect
four elements of ${\cal F}_{i_1i_2j_1j_2k_1k_2}$ into a tree.
In this case, the four elements of ${\cal F}_{i_1i_2j_1j_2k_1k_2}$ are
different ${\cal B}_1$-equivalence classes, which are obviously
singletons. Therefore the set of indispensable monomials for this
problem is
\[
\{u_{i_1j_1k_1}u_{i_2j_2k_2},\ u_{i_1j_1k_2}u_{i_2j_2k_1},\
u_{i_1j_2k_1}u_{i_2j_1k_2},\ u_{i_1j_2k_2}u_{i_2j_1k_1},\ i_1\neq
i_2,j_1\neq j_2,k_1\neq k_2\}
\]
in addition to the positive and negative components of ${\cal
B}_{{\rm IDP}}$. 
Figure \ref{fig:4-element-monomials} illustrates the fiber ${\cal
F}_{i_1i_2j_1j_2k_1k_2}$. 

\begin{figure}[htbp]
\begin{center}
\includegraphics[width=50mm]{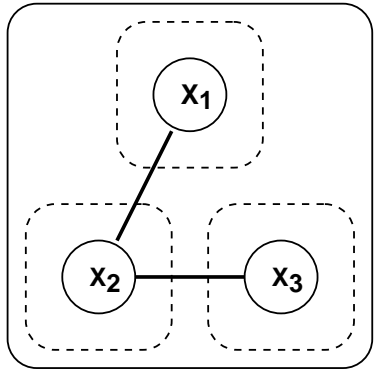}
\caption{Illustration of the 3-element fiber ${\cal
F}_{i_1i_2i_3i_4}$ of the Hardy-Weinberg models. The three circles are three
 monomials, 
 $\Bx_1 = u_{i_1i_2}u_{i_3i_4}, \Bx_2 = u_{i_1i_3}u_{i_2i_4},
\Bx_3 = u_{i_1i_4}u_{i_2i_3}$. 
 The thick
 lines are the two dispensable moves, $\{\Bx_1-\Bx_2,
 \Bx_2-\Bx_3\}$.}
\label{fig:3-element-monomials}
\end{center}
\end{figure}

\paragraph*{Hardy-Weinberg model.}\ Another example considered in
\cite{takemura-aoki-2004aism} is the
Hardy-Weinberg model for $I$ alleles, i.e.,
\[
\Bx =
(x_{11},x_{12},\ldots,x_{1I},x_{22},x_{23},\ldots,x_{2I},x_{33},\ldots,x_{II})'
\]
and
\[
\matrixA = (\matrixA_I\ \matrixA_{I-1}\ \cdots\ \matrixA_1),\quad
\matrixA_k = \left(O_{k\times(I-k)}\
  B_k'\right)',
\]
where $B_k$ is the following $k\times k$ square matrix
\[
B_k = \left[\begin{array}{ccccc}
    2      & 1      & 1 & \cdots & 1\\
    0      & 1      & 0 & \cdots & 0\\
    0      & 0      & 1 &        & 0\\
    \vdots & \vdots &   & \ddots & \vdots\\
    0      & 0      & \cdots & 0 & 1
\end{array}\right].
\]
As is stated in \cite{takemura-aoki-2004aism}, a minimal Markov basis for
this case is not unique, and the minimum fiber Markov basis is:
\[
\begin{array}{l}
{\cal B}_{\rm{MF}} = {\cal B}_{\rm{IDP}} \cup {\cal B}^{*},\\
{\cal B}_{\rm{IDP}} = \{u_{i_1i_1}u_{i_2i_3} - u_{i_1i_2}u_{i_1i_3},\
u_{i_1i_1}u_{i_2i_2} - u_{i_1i_2}^2\},\\
{\cal B}^{*} = \{u_{i_1i_2}u_{i_3i_4} - u_{i_1i_3}u_{i_2i_4},\
 u_{i_1i_2}u_{i_3i_4} -
u_{i_1i_4}u_{i_2i_3},\ u_{i_1i_3}u_{i_2i_4} - u_{i_1i_4}u_{i_2i_3}\},
\end{array}
\]
where $i_1,i_2,i_3,i_4$ are all distinct, and $u_{ij} = u_{ji}$ for $i
> j$. Here, 
${\cal B}^{*}$ is the set of all degree 2 moves which connect all of the
  elements of the three-element fiber
\[
{\cal F}_{i_1i_2i_3i_4} = \{u_{i_1i_2}u_{i_3i_4},\
u_{i_1i_3}u_{i_2i_4},\ u_{i_1i_4}u_{i_2i_3}\}.
\]
Again, these three elements of ${\cal F}_{i_1i_2i_3i_4}$ form
singleton ${\cal B}_1$-equivalence classes of it by themselves, and are
indispensable monomials.
Figure \ref{fig:3-element-monomials} illustrates the fiber ${\cal
F}_{i_1i_2i_3i_4}$. 

\subsection{Examples of Case 3}
Some examples for this case are found in the hierarchical models of $2\times
2\times 2\times 2$ contingency tables considered in 
\cite{aoki-takemura-2003metr}. First we will show one of them
as an example of Case 3. By modifying the example, 
we will show another example of the situation
considered in Proposition \ref{prop:3}, i.e., the situation that some
dispensable moves contain
both indispensable and dispensable monomials as their positive and
negative parts.

\begin{figure}[htbp]
\begin{center}
\includegraphics[width=50mm]{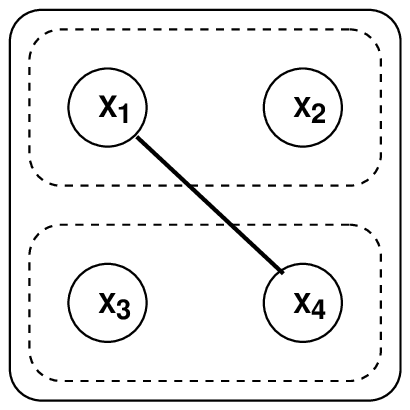}
\caption{Illustration of the 4-element fiber ${\cal
F}^1$ of $12/13/23/34$ model of $2\times 2\times 2\times 2$ contingency
 tables. The four
 monomials are
$\Bx_1 = u_{1111}u_{1221}u_{2121}u_{2212}, \Bx_2 = 
 u_{1112}u_{1221}u_{2121}u_{2211}, \Bx_3 =
 u_{1121}u_{1211}u_{2112}u_{2221}, \mbox{ and }\Bx_4 = 
u_{1121}u_{1212}u_{2111}u_{2221}$. $\{\Bx_1,\Bx_2\}$ and $\{\Bx_3,\Bx_4\}$
 form
 ${\cal B}_{3}$-equivalence classes of the fiber. The thick
 line is a dispensable move, $\Bx_1-\Bx_4$.}
\label{fig:4-element-monomials-2x2x2x2}
\end{center}
\end{figure}

\begin{figure}[htbp]
\begin{center}
\includegraphics[width=80mm]{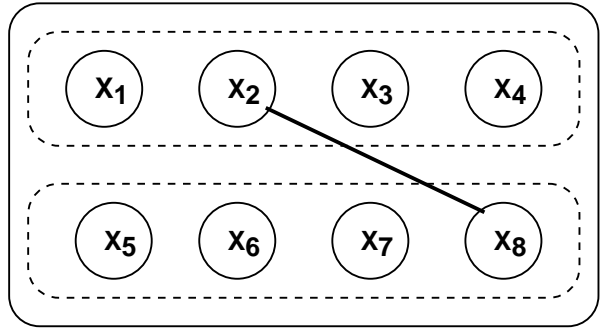}
\caption{Illustration of the 8-element fiber ${\cal
F}^2$ of $12/13/23/34$ model of $2\times 2\times 2\times 2$ contingency
 tables. The eight
 monomials are $\Bx_1 = u_{1111}u_{1221}u_{2122}u_{2212}$, $\Bx_2 =
 u_{1112}u_{1222}u_{2121}u_{2211}$, $\Bx_3 = u_{1111}u_{1222}u_{2121}u_{2212}$,
$\Bx_4 = u_{1112}u_{1221}u_{2122}u_{2211}$, $\Bx_5 =
 u_{1121}u_{1211}u_{2112}u_{2222}$, $\Bx_6 =
 u_{1122}u_{1212}u_{2111}u_{2221}$,
$\Bx_7 = u_{1121}u_{1212}u_{2111}u_{2222}$,  and $\Bx_8 =
 u_{1122}u_{1211}u_{2112}u_{2221}$. 
$\{\Bx_1,\Bx_2,\Bx_3,\Bx_4\}$
 and $\{\Bx_5,\Bx_6,\Bx_7,\Bx_8\}$
 form
 ${\cal B}_{3}$-equivalence classes of the fiber. The thick
 line is a dispensable move, $\Bx_2-\Bx_8$.}
\label{fig:8-element-monomials}
\end{center}
\end{figure}

\paragraph*{$12/13/23/34$ model of $2\times 2\times 2\times 2$
    contingency tables.}\ 
Let $\Bx$ be a frequency vector for $2\times
    2\times 2\times 2$ contingency tables ($p = 16$). We write
    indeterminates with respect to a lexicographic order as
\[
\Bu = \{u_{1111},u_{1112},u_{1121},u_{1122},u_{1211},\ldots,u_{2222}\}.
\]
Consider the model of $d = 9$ given as
\[
 \matrixA = \left[\begin{array}{cccccccccccccccc}
1 & 1 & 1 & 1 & 1 & 1 & 1 & 1 & 1 & 1 & 1 & 1 & 1 & 1 & 1 & 1\\
1 & 1 & 1 & 1 & 1 & 1 & 1 & 1 & 0 & 0 & 0 & 0 & 0 & 0 & 0 & 0\\
1 & 1 & 1 & 1 & 0 & 0 & 0 & 0 & 1 & 1 & 1 & 1 & 0 & 0 & 0 & 0\\
1 & 1 & 0 & 0 & 1 & 1 & 0 & 0 & 1 & 1 & 0 & 0 & 1 & 1 & 0 & 0\\
1 & 0 & 1 & 0 & 1 & 0 & 1 & 0 & 1 & 0 & 1 & 0 & 1 & 0 & 1 & 0\\
1 & 1 & 1 & 1 & 0 & 0 & 0 & 0 & 0 & 0 & 0 & 0 & 0 & 0 & 0 & 0\\
1 & 1 & 0 & 0 & 1 & 1 & 0 & 0 & 0 & 0 & 0 & 0 & 0 & 0 & 0 & 0\\
1 & 1 & 0 & 0 & 0 & 0 & 0 & 0 & 1 & 1 & 0 & 0 & 0 & 0 & 0 & 0\\
1 & 0 & 0 & 0 & 1 & 0 & 0 & 0 & 1 & 0 & 0 & 0 & 1 & 0 & 0 & 0
\end{array}
\right].
\]
For this case, there are 12 indispensable moves of degree 2 and 4
indispensable moves of degree 4, but the set of indispensable moves does
not form a Markov basis. In addition to the indispensable moves,
we have to consider moves connecting 4-element fiber
\[
 {\cal F}^1 = \{u_{1111}u_{1221}u_{2121}u_{2212},
 u_{1112}u_{1221}u_{2121}u_{2211},u_{1121}u_{1211}u_{2112}u_{2221},
 u_{1121}u_{1212}u_{2111}u_{2221}\}
\]
and 8-element fiber
\[\begin{array}{l}
 {\cal F}^2 = \{u_{1111}u_{1221}u_{2122}u_{2212},
 u_{1112}u_{1222}u_{2121}u_{2211}, u_{1111}u_{1222}u_{2121}u_{2212},
 u_{1112}u_{1221}u_{2122}u_{2211},\\
\hspace*{15mm} u_{1121}u_{1211}u_{2112}u_{2222},
 u_{1122}u_{1212}u_{2111}u_{2221},
 u_{1121}u_{1212}u_{2111}u_{2222}, u_{1122}u_{1211}u_{2112}u_{2221} \}.
\end{array}
\]
For ${\cal F}^1$, we have ${\cal B}_3$-equivalence classes of it as
\[
 {\cal F}^1 = \{u_{1111}u_{1221}u_{2121}u_{2212},
 u_{1112}u_{1221}u_{2121}u_{2211}\}\cup\{
u_{1121}u_{1211}u_{2112}u_{2221},
 u_{1121}u_{1212}u_{2111}u_{2221}\}.
\]
Therefore these 4 elements are dispensable monomials. In fact, we can
find an element of minimal Markov basis not containing
$u_{1111}u_{1221}u_{2121}u_{2212}$, for example, as
\[
 u_{1112}u_{1221}u_{2121}u_{2211} - 
u_{1121}u_{1211}u_{2112}u_{2221}.
\]
Similarly for ${\cal F}^2$, ${\cal B}_3$-equivalence classes of it are
given as
\[
 \begin{array}{l}
 {\cal F}^2 = \{u_{1111}u_{1221}u_{2122}u_{2212},
 u_{1112}u_{1222}u_{2121}u_{2211}, u_{1111}u_{1222}u_{2121}u_{2212},
 u_{1112}u_{1221}u_{2122}u_{2211}\},\\
\hspace*{15mm} \cup \{u_{1121}u_{1211}u_{2112}u_{2222},
 u_{1122}u_{1212}u_{2111}u_{2221},
 u_{1121}u_{1212}u_{2111}u_{2222}, u_{1122}u_{1211}u_{2112}u_{2221} \}.
\end{array}
\]
Figure \ref{fig:4-element-monomials-2x2x2x2} and Figure
\ref{fig:8-element-monomials} illustrate the fiber ${\cal F}^1$ and
${\cal F}^2$, respectively. 


\begin{figure}[htbp]
\begin{center}
\includegraphics[width=50mm]{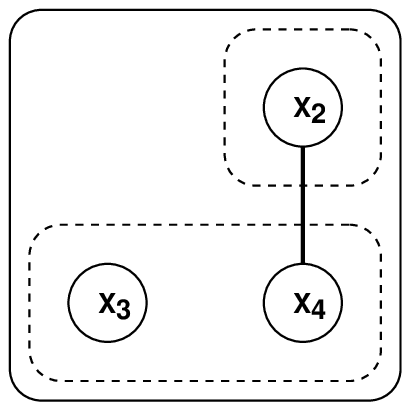}
\caption{Illustration of the 3-element fiber ${\cal
F}^{1*}$ of $12/13/23/34$ model of $2\times 2\times 2\times 2$ contingency
 tables with structural zero cell. This fiber is constructed by removing
 element $\Bx_1$ from ${\cal F}^1$. In this case, 
$\{\Bx_2\}$ and $\{\Bx_3,\Bx_4\}$ form
 ${\cal B}_{3}$-equivalence classes of the fiber. The thick
 line means a choice of a dispensable move, $\Bx_2-\Bx_4$. Another
 possibility of constructing a minimal Markov basis is to choose a
 dispensable move, $\Bx_2 - \Bx_3$. The monomial $\Bx_2$ is
 included in any minimal Markov basis, and is an indispensable monomial.}
\label{fig:3-element-monomials-2x2x2x2}
\end{center}
\end{figure}

\paragraph*{$12/13/23/34$ model of $2\times 2\times 2\times 2$
    contingency tables with a structural zero cell.}\ We modify the
    previous example by introducing
    a structural zero cell, $x_{111} \equiv 0$. This situation
    corresponds to removing the indeterminate $\Bu_{1111}$ and the first
    column of $\matrixA$ as
\[
\Bu = \{u_{1112},u_{1121},u_{1122},u_{1211},\ldots,u_{2222}\},
\]
\[
 \matrixA = \left[\begin{array}{ccccccccccccccc}
1 & 1 & 1 & 1 & 1 & 1 & 1 & 1 & 1 & 1 & 1 & 1 & 1 & 1 & 1\\
1 & 1 & 1 & 1 & 1 & 1 & 1 & 0 & 0 & 0 & 0 & 0 & 0 & 0 & 0\\
1 & 1 & 1 & 0 & 0 & 0 & 0 & 1 & 1 & 1 & 1 & 0 & 0 & 0 & 0\\
1 & 0 & 0 & 1 & 1 & 0 & 0 & 1 & 1 & 0 & 0 & 1 & 1 & 0 & 0\\
0 & 1 & 0 & 1 & 0 & 1 & 0 & 1 & 0 & 1 & 0 & 1 & 0 & 1 & 0\\
1 & 1 & 1 & 0 & 0 & 0 & 0 & 0 & 0 & 0 & 0 & 0 & 0 & 0 & 0\\
1 & 0 & 0 & 1 & 1 & 0 & 0 & 0 & 0 & 0 & 0 & 0 & 0 & 0 & 0\\
1 & 0 & 0 & 0 & 0 & 0 & 0 & 1 & 1 & 0 & 0 & 0 & 0 & 0 & 0\\
0 & 0 & 0 & 1 & 0 & 0 & 0 & 1 & 0 & 0 & 0 & 1 & 0 & 0 & 0\\
\end{array}
\right].
\]
In this case, the fiber ${\cal F}^1$ in the previous example is modified
to a 3-element fiber, 
\[
 {\cal F}^{1*} = \{
 u_{1112}u_{1221}u_{2121}u_{2211},u_{1121}u_{1211}u_{2112}u_{2221},
 u_{1121}u_{1212}u_{2111}u_{2221}\}.
\]
Since ${\cal F}^{1*}$ has still different ${\cal B}_3$-equivalence
classes, we have to consider moves connecting the elements of ${\cal
F}^{1*}$ to construct a minimal Markov basis. In this case, ${\cal
B}_3$-equivalence classes of ${\cal F}^{1*}$ are given as
\[
 {\cal F}^{1*} = \{
u_{1112}u_{1221}u_{2121}u_{2211}\}\cup\{
u_{1121}u_{1211}u_{2112}u_{2221},
 u_{1121}u_{1212}u_{2111}u_{2221}\}.
\]
Therefore a minimal Markov basis for this problem has to contain either 
\[
u_{1112}u_{1221}u_{2121}u_{2211} - u_{1121}u_{1211}u_{2112}u_{2221}
\]
or
\[
u_{1112}u_{1221}u_{2121}u_{2211} - u_{1121}u_{1212}u_{2111}u_{2221}.
\]
The above two moves are dispensable, and the negative parts of both
moves are also dispensable monomials, whereas the positive part,
$u_{1112}u_{1221}u_{2121}u_{2211}$, is an indispensable monomial.
Figure \ref{fig:3-element-monomials-2x2x2x2} 
illustrates the fiber ${\cal F}^{1*}$.

\section{Some discussions}
\label{sec:discussions}
In this paper, the concept of indispensable monomials is introduced, by
extending the notion of indispensable binomials. Both in the framework of
Markov bases and toric ideals, the indispensable monomial plays an
important role since it has to be included in all Markov
bases or generators of toric ideals. 
It is true that enumerating indispensable monomials is as difficult
as enumerating indispensable binomials.

Note that, by the notion of indispensable monomials, we can 
characterize a dispensable binomial as (i) a
difference of two dispensable monomials, (ii) a difference of
dispensable and indispensable monomials, or (iii) a difference of two
indispensable monomials.  The situations where each case arises are
shown in Proposition \ref{prop:3} and in Proposition \ref{prop:1}.  We
have found some examples for the case (ii) by introducing some
structural zero cells for the case (i),

The enumeration of indispensable monomials seems very important problem, 
since it can lead directly to the enumeration of indispensable
binomials. In addition, it also gives the fibers of the special structure
that it contains at least one singleton equivalence class. Moreover, by
finding dispensable binomials which are differences of 
two indispensable monomials, we can find all fibers that only contain
singleton equivalence classes.


\begin{thebibliography}{10}

\bibitem{aoki-takemura-2003metr}
Satoshi Aoki and Akimichi Takemura.
\newblock Invariant minimal markov basis for sampling contingency tables with
  fixed marginals.
\newblock METR 2003-25, 2003.
\newblock Submitted for publication.

\bibitem{aoki-takemura-2003anz}
Satoshi Aoki and Akimichi Takemura.
\newblock Minimal basis for a connected {M}arkov chain over {$3\times 3\times
  K$} contingency tables with fixed two-dimensional marginals.
\newblock {\em Aust. N. Z. J. Stat.}, 45(2):229--249, 2003.

\bibitem{diaconis-sturmfels}
Persi Diaconis and Bernd Sturmfels.
\newblock Algebraic algorithms for sampling from conditional distributions.
\newblock {\em Ann. Statist.}, 26(1):363--397, 1998.

\bibitem{Hosten2004}
Serkan Ho{\c{s}}ten, Amit Khetan, and Bernd Sturmfels.
\newblock Solving the likelihood equations.
\newblock {\em Foundations of Computational Mathematics}, 2004.
\newblock To appear.

\bibitem{Lehmann-tsh-3rd}
E.~L. Lehmann and Joseph~P. Romano.
\newblock {\em Testing statistical hypotheses}.
\newblock Springer Texts in Statistics. Springer, New York, third edition,
  2005.

\bibitem{miller-sturmfels}
Ezra Miller and Bernd Sturmfels.
\newblock {\em Combinatorial Commutative Algebra}, volume 227 of {\em Graduate
  Texts in Mathematics}.
\newblock Springer-Verlag, New York, 2005.

\bibitem{ohsugi-hibi-indispensable}
Hidefumi Ohsugi and Takayuki Hibi.
\newblock Indispensable binomials of finite graphs.
\newblock {\em Journal of Algebra and Its Applications}, 4(4):421--434, 2005.

\bibitem{ohsugi-hibi-contingency-tables-2005}
Hidefumi Ohsugi and Takayuki Hibi.
\newblock Toric ideals arising from contingency tables.
\newblock 2005.
\newblock Proceedings of the Ramanujan Mathematical Society's Lecture Notes
  Series, to appear.

\bibitem{Pachter2004a}
Lior Pachter and Bernd Sturmfels.
\newblock Parametric inference for biological sequence analysis.
\newblock {\em Proc Natl Acad Sci U S A}, 101(46):16138--43, 2004.

\bibitem{ascb}
Lior Pachter and Bernd Sturmfels.
\newblock {\em Algebraic Statistics for Computational Biology}.
\newblock Cambridge University Press, Cambridge, UK, 2005.

\bibitem{Pistone2001}
Giovanni Pistone, Eva Riccomagno, and Henry~P. Wynn.
\newblock {\em Algebraic Statistics: Computational Commutative Algebra in
  Statistics}.
\newblock Chapman \& Hall Ltd, Boca Raton, 2001.

\bibitem{sturmfels1996}
Bernd Sturmfels.
\newblock {\em Gr\"obner Bases and Convex Polytopes}, volume~8 of {\em
  University Lecture Series}.
\newblock American Mathematical Society, Providence, RI, 1996.

\bibitem{Sullivant2005}
Seth Sullivant.
\newblock {Compressed polytopes and statistical disclosure limitation}.
\newblock arXiv:math.CO/0412535.

\bibitem{takemura-aoki-2004aism}
Akimichi Takemura and Satoshi Aoki.
\newblock Some characterizations of minimal {M}arkov basis for sampling from
  discrete conditional distributions.
\newblock {\em Ann. Inst. Statist. Math.}, 56(1):1--17, 2004.

\bibitem{takemura-aoki-2005bernoulli}
Akimichi Takemura and Satoshi Aoki.
\newblock Distance reducing markov bases for sampling from a discrete sample
  space.
\newblock {\em Bernoulli}, 11(5):793--813, 2005.

\end{thebibliography}
\end{document}